\documentclass[12pt,a4paper]{article}
\usepackage{amssymb,amsmath,amsthm}
\newtheorem{theorem}{\indent {\sc Theorem}}
\begin{document}

{\large

\begin{center}
\textbf{\Large{Ternary Goldbach's Problem Involving Primes\\ of a
Special Type}}
\end{center}

\bigskip

\begin{center}
\textbf{\large{Sergey A. Gritsenko and Natalya N. Motkina}}
\end{center}
\bigskip
\begin{abstract}
Let $\eta$ be a quadratic irrationality. The variant of a ternary
problem of Goldbach involving primes such that $a<\{\eta p\}<b$,
where $a$ and $b$ are arbitrary numbers of the interval $(0,1)$,
solved in this paper.
\end{abstract}
\bigskip

\textbf{1. Introduction} Suppose that $k\ge 2$ and $n\ge 1$ are
naturals. Consider the equation
\begin{equation}\label{e1}
p_1^n+p_2^n+\cdots +p_k^n=N
\end{equation}
in primes $p_1,\ p_2,\ldots ,p_k$.

Several classical problems in number theory can be reduced to the
question on the number of solutions of the equation  (\ref{e1}).
For example if $n=1$ and $k=2$ or $3$ then (\ref{e1}) is the
equation of Goldbach; if $n\ge 3$ then (\ref{e1}) is the equation
of Waring--Goldbach.

Let $\mathcal{P}$ be a subset of the set of primes,
$\mu=\lim\limits_{n\to\infty}\pi^{-1}(N)\sum\limits_{\substack{p\le
N\\p\in \mathcal{P}}}1$ be the "density" of $\mathcal{P}$,
$0<\mu<1$.

Let $J_{k,n}(N)$ be the number of solutions of (\ref{e1}) in
primes of the set $\mathcal{P}$, and $I_{k,n}(N)$ be the number of
solutions of (\ref{e1}) in arbitrary primes.

Sometimes the equality
\begin{equation}\label{e2}
J_{k,n}(N)\sim\mu^k I_{k,n}(N)
\end{equation}
holds.

For example, if $\mathcal{P}=\{p\ |\ a<\{\sqrt{p}\}< b \}$,
$0<a<b<1$ then provided $n=1$, $k=3$, and also provided $n\ge 3$,
$k\gg n^2\log n$ the equality (\ref{e2}) holds (s. \cite {G}).

In present paper we solve the additive problem for which the
equality (\ref{e2}) does not hold.

Let $\eta$ be a quadratic irrationality, $0<a<b<1$,
$\mathcal{P}=\{p\ |\ a<\{\eta p\}<b\}$, $n=1$. Then $\mu=b-a$ (s.
\cite {V}).

The main result of the present paper is the following theorem
which
 shows that the equality  $J_{3,1}(N)\sim (b-a)^3 I_{3,1}(N)$ does not hold.

Note that (s. [1])
$$ I_{3,1}(N)=\sigma (N) \frac{N^2}{2(\log
N)^3}+O\left(\frac{N^2}{(\log N)^4}\right),
$$
where
$$
\sigma(N)=\prod_{p}\left(1+\frac{1}{(p-1)^3}\right)\prod_{p\backslash
N}\left(1+\frac{1}{p^2-3p+3}\right).
$$

\begin{theorem}
Suppose
$$\sigma(N,a,b)=\sum_{|m|<\infty} e^{2\pi i m(\eta
N-1.5(a+b))} \frac{\sin^3 \pi m (b-a)}{\pi ^3 m^3}.
$$
Then, for any $C>0$, formula
$$
J_{3,1}(N)=I_{3,1}(N)\sigma(N,a,b)+O(N^2 \ln^{-C}N)
$$
holds.
\end{theorem}

{\textbf{Outline of proof.} Define the function
$$\psi_0 (x)= \bigl\{
\begin{array}{ll} 1,\quad \emph{if}\quad & a<x<b,\\ 0,
\quad \emph{if}\quad & 0\leq x\leq a \quad \emph{or} \quad b\leq x \leq 1
\end{array} $$
and extend it along number axis with period 1. Suppose
$$
S_0(x)=\sum_{p\leq N}\psi_0 (\eta p) e^{2 \pi i x p};
$$
then
$$
J_{3,1}(N)=\int_0^1 S^3_0(x)e^{-2 \pi i x N}dx.
$$
In Vinogradov's  lemma on ''containers''(s. \cite {V}, p. 22), we
select $r=[\ln N]$, $\Delta = \ln^{-A_1}N$, where $A_1>0$ is
sufficiently large constant. Denote the function $\psi$ from the
lemma on the ''containers'' by $\psi_1$ with $\alpha=a+\Delta/2$
and $\beta=b-\Delta/2$. Set in the lemma $\alpha=a-\Delta/2$,
$\beta=b+\Delta/2$ and denote the corresponding function $\psi$ by
$\psi_2$. Define
\begin{equation}\label{1}
J_k(N)=\int_0^1 \left(\sum_{p\leq N}\psi_k (\eta p) e^{2 \pi i x
p}\right)^3e^{-2 \pi i x N}dx, \quad k=1,2.
\end{equation}
From the properties of $\psi_1(x)$ and $\psi_2(x)$ there follows
\begin{equation}\label{2}
J_1(N)\leq J_{3,1}(N) \leq J_2(N).
\end{equation}
We shall get an asymptotic formulas for $J_1(N)$ and $J_2(N)$ with
identical main terms. In the above integral (\ref{1}) we replace
$\psi_k (\eta p)$ by the Fourier series
$$
\psi_k (\eta p)=\sum_{|m|\leq r\Delta ^ {-1}}c_k(m)e^{2 \pi i m
\eta p}+O(N^{-\ln \pi});
$$
thus
$$
J_k(N)=\int_0^1 \left(\sum_{|m|\leq r\Delta ^
{-1}}c_k(m)\sum_{p\leq N} e^{2 \pi i (x+m\eta ) p} \right)^3e^{-2
\pi i x N}dx +O(N^{2-\ln \pi}) $$ $$=\sum_{|m_1|\leq r\Delta ^
{-1}} c_k(m_1) \sum_{|m_2|\leq r\Delta ^ {-1}}c_k(m_2)
\sum_{|m_3|\leq r\Delta ^ {-1}} c_k(m_3) \int_0^1 \sum_{p_1\leq N}
e^{2\pi i(x+ m_1\eta )p_1}\cdot $$ $$\cdot \sum_{p_2\leq N}e^{2\pi
i(x+ m_2\eta )p_2} \sum_{p_3\leq N} e^{2\pi i(x+
m_3\eta)p_3}e^{-2\pi ix N} dx +O(N^{2-\ln \pi}).
$$

Note that
$$
\sum_{|m|\leq r\Delta ^ {-1}}c^3_k(m)e^{2\pi in\eta
N}\sum_{p_1\leq N} \sum_{p_2\leq N} \sum_{p_3\leq N} \int_0^1
e^{2\pi i(x+m\eta) (p_1+p_2+p_3-N)}dx=
$$
$$=I_{3,1}(N)(\sigma(N,a,b)+O({\Delta})).$$

Consider the sets such that $(m_1,m_2,m_3)\neq (m,m,m)$  and the
integrals
$$
I(N,m_1,m_2,m_3)=\int_0^1 |S(x+ m_1\eta)||S(x+ m_2\eta)||S(x+
m_3\eta)|dx,
$$
where
$$
S(x)=\sum_{p\leq N} e^{2\pi ixp}.
$$
Without loss of generality it can be assumed that $m_1<m_2$. By
$t$ denote  $x+ m_1 \eta$. Since the integrand has a period $1$,
we can take $t$ to lie in the range $E=[-1/\tau;1-1/\tau)$, where
$\tau=N\ln^{-B}N$, $B>2C+8$. Suppose  $d,q\in \mathbb{Z}$ are such
that
$$
t=\frac{d}{q}+\frac{\theta_1}{q\tau},\quad (d,q)=1,\quad 1\leq
q\leq \tau,\quad |\theta_1|<1.
$$
Divide the interval $E$ into two intervals $E_1$ and $E_2$ for
which $q\leq \ln^A N$ and $\ln^A N <q\leq \tau$ respectively with
a fixed $A>2C+8$. On writing $m=m_2-m_1$, $m'=m_3-m_1$, we can
express the above integral as
$$
\int_{-1/\tau}^{1-1/\tau} |S(t)||S(t+m\eta)||S(t+m' \eta)|dt$$ and
get the estimate
\begin{equation}\label{3}
\ll \pi(N)(\max_{t\in E_1}|S(t+m\eta)|+\max_{t\in E_2}|S(t|).
\end{equation}
It is well known (s. \cite {Vh}, p. 35) that
\begin{equation}\label{4}
\max_{t\in E_2}|S(t)|\ll N^{1/2}{\tau}^{1/2}\ln^4 N.
\end{equation}
Suppose  $D,Q\in \mathbb{Z}$ are such that
$$
\eta=\frac{D}{Q}+\frac{\theta_2}{Q^2},\quad (D,Q)=1,\quad 1\leq
Q\leq \tau_1,\quad |\theta_2|<1.
$$
From $\eta$ is the quadratic irrationality it follows that $\tau
_1\geq Q>c(\alpha)\tau _1$. Choosing $\tau_1=\sqrt{ \tau}$, we can
get the estimate
\begin{equation}\label{5}
\max_{t\in E_1}|S(t+m\eta )|\ll N^{1/2}{\tau}^{1/2}\ln^4 N.
\end{equation}
Now the theorem  follows directly from formulas  (\ref{2}) ---
(\ref{5})


\begin{thebibliography}{3}

\bibitem{V}
I.M.Vinogradov, \textit{The method of trigonometrical Sums} [in
Rusian], Nauka, Moscow (1983)

\bibitem{G}
S.A. Gritsenko, \textit{Three additive problems}, Izvestiya
Mathematics, 1993, 41 (3), 447--464

\bibitem{Vh}
R.C.Vaughan, \textit{The Hardy-Littlewood Method} [in Russian],
Mir, Moscow (1985)
\end{thebibliography}
\end{document}